\documentclass[10pt]{article}

\usepackage{amsmath,amscd,amsfonts,amsthm}

\newcommand{\shrinkmargins}[1]{
  \addtolength{\textheight}{#1\topmargin}
  \addtolength{\textheight}{#1\topmargin}
  \addtolength{\textwidth}{#1\oddsidemargin}
  \addtolength{\textwidth}{#1\evensidemargin}
  \addtolength{\topmargin}{-#1\topmargin}
  \addtolength{\oddsidemargin}{-#1\oddsidemargin}
  \addtolength{\evensidemargin}{-#1\evensidemargin}
  }

\shrinkmargins{.7}

\DeclareMathOperator{\Hom}{Hom}

\DeclareMathOperator{\PGL}{PGL}
\DeclareMathOperator{\GL}{GL}

\DeclareMathOperator{\Spec}{Spec}

\DeclareMathOperator{\Gal}{Gal}

\DeclareMathOperator{\Res}{Res}

\DeclareMathOperator{\Jac}{Jac}
\DeclareMathOperator{\Cot}{Cot}

\newcommand{\field}[1]{\mathbb{#1}}

\newcommand{\Q}{\field{Q}}
\newcommand{\Qbar}{\bar{\Q}}

\newcommand{\Z}{\field{Z}}
\newcommand{\F}{\field{F}}
\newcommand{\Fp}{\field{F}_p}

\newcommand{\C}{\field{C}}
\renewcommand{\P}{\field{P}}
\newcommand{\ord}{\mbox{ord}}

\newcommand{\CC} {\mathcal{C}}

\newcommand{\ra}{\rightarrow}

\newcommand{\OO}{\mathcal{O}}

\newcommand{\ic}[1]{\mathfrak{#1}}

\newcommand{\mat}[4]{\left[\begin{array}{cc}#1 & #2 \\
                                         #3 & #4\end{array}\right]}

\newcommand{\tensor} {\otimes}

\newcommand{\GG}{\mathcal{G}}

\newcommand{\GalK}{\Gal(\bar{K}/K)}
\newcommand{\GalQ}{\Gal(\bar{\Q}/\Q)}
\newcommand{\set}[1]{\{#1\}}

\newcommand{\FF}{\mathcal{F}}
\newcommand{\notdiv}{\not | \,}
\newcommand{\beq}{\begin{displaymath}}
\newcommand{\eeq}{\end{displaymath}}
\newcommand{\beqn}{\begin{equation}}
\newcommand{\eeqn}{\end{equation}}

\theoremstyle{plain}
\newtheorem{thm}{Theorem}[section]
\newtheorem{prop}[thm]{Proposition}

\newtheorem{lem}[thm]{Lemma}
\newtheorem*{intro}{Theorem}

\theoremstyle{definition}
\newtheorem{defn}[thm]{Definition}

\theoremstyle{remark}
\newtheorem{rem}[thm]{Remark}

\title{Galois representations attached to $\Q$-curves and the
generalized Fermat equation $A^4 + B^2 = C^p$}

\author{Jordan S. Ellenberg \footnote{Partially supported by NSA Young Investigator Grant
  MDA905-02-1-0097.} \\ Princeton University \\ \texttt{ellenber@math.princeton.edu}}
\date{22 Jul 2003}

\begin{document}

\maketitle

\begin{abstract}  We prove that the equation $A^4 + B^2 = C^p$ has no
solutions in coprime positive integers when $p \geq 211$.  The main
step is to show that, for all sufficiently large primes $p$, every
$\Q$-curve over an imaginary quadratic field $K$ with a prime of bad
reduction greater than $6$ has a surjective mod $p$ Galois
representation.  The bound on $p$ depends on $K$ and the degree of the isogeny
between $E$ and its Galois conjugate, but is independent of the choice
of $E$.  The proof of this theorem combines geometric arguments due to 
Mazur, Momose, Darmon, and Merel with an analytic estimate of the
average special values of certain $L$-functions.
\end{abstract}

\section{Introduction}

The resolution of the Fermat problem has demonstrated a close
relationship between the solutions of Diophantine equations and the
arithmetic of abelian varieties over number fields.  It remains far
from clear which Diophantine equations can be productively studied
along the lines developed by Frey, Hellegouarch, Serre, Ribet, Wiles,
and Taylor.

In particular, one wonders whether modular abelian varieties can
address the classical problem of describing all solutions to the {\em
generalized Fermat equation} 
\begin{equation}
A^p + B^q = C^r
\label{eq:gfe}
\end{equation}

in coprime integers $A,B,C$.  Darmon and Granville~\cite{darm:dagr}
have proved that \eqref{eq:gfe} has only finitely many solutions for
any particular $p,q,r$ satisfying $1/p + 1/q + 1/r < 1$.  It is
conjectured that \eqref{eq:gfe} has only finitely many solutions,
excepting $1^p + 2^3 = 3^2$, even if $p,q,r$ are allowed to vary
(still subject to $1/p + 1/q + 1/r < 1$.)

Certain special cases of \eqref{eq:gfe} and similar equations have
been treated by Darmon, Merel, and
Ribet(\cite{darm:imrn},\cite{darm:wind},\cite{ribe:ap2bpcp}) using
elliptic curves over $\Q$, and by Bruin~\cite{brui:gf} using Chabauty
methods.  In ~\cite{darm:rigid}, Darmon discusses
the relationship between more general cases of \eqref{eq:gfe} and
as-yet-unproved conjectures about the Galois representations attached
to Hilbert-Blumenthal abelian varieties over number fields.

Our goal in the present paper is twofold.  Our main motivation (or, as
Darmon and Merel put it in~\cite{darm:wind}, our ``excuse'') is to
prove the following Diophantine theorem:

\begin{intro}[Theorem $4.1$]
Suppose $A,B,C$ are coprime integers such that 
\begin{equation}
A^4 + B^2 = C^p
\label{eq:intro}
\end{equation}
and $p \geq 211$.  Then $AB = 0$.
\end{intro}

We will attach an elliptic curve $E$ to any solution to \eqref{eq:intro}.
However, $E$ will be defined not over $\Q$, but over
$\Q[i]$, and it will be isogenous to its Galois conjugate.  An
elliptic curve, like $E$, whose isogeny class is defined over $\Q$ is
called a {\em $\Q$-curve}.  (The idea of studying \eqref{eq:intro} by
means of $\Q$-curves was arrived at independently by Darmon in
\cite{darm:sc}.)

In order to prove Theorem 4.1, it is necessary 
to bring our knowledge of the arithmetic of $\Q$-curves more in line
with our knowledge about elliptic curves defined over $\Q$.  Our
second goal in this paper is to prove a theorem on surjectivity of mod
$p$ representations attached to $\Q$-curves with non-integral
$j$-invariant.  If $K$ is a quadratic
number field, a $\Q$-curve $E/K$ of degree $d$ is an elliptic curve
over $K$ which admits a cyclic isogeny of degree $d$ to its Galois
conjugate.

\begin{intro}[Theorem 3.14]
Let $K$ be an imaginary quadratic field and $d$ a square-free positive integer.
There exists an effective constant $M_{K,d}$ such that, for all primes $p >
M_{K,d}$ and all $\Q$-curves $E/K$ of degree $d$, either 
\begin{itemize}
\item
the representation
\beq
\P\bar{\rho}_{E,p}: \GalQ \ra \PGL_2(\F_p)
\eeq
is surjective, or
\item $E$ has potentially good reduction at all primes not dividing
$6$.
\end{itemize}
\end{intro}

Theorem 3.14 plays the role that Mazur's theorem~\cite{mazu:ripd} does
in the solution of Fermat's problem.  It is worth remarking that the
analogue of Theorem 3.14 for elliptic curves over $\Q$ is still conjectural.

We also need a modularity
theorem for $\Q$-curves; we have proved the result we need in an
earlier paper with C. Skinner~\cite{elle:qcurves}.

Much of the proof of Theorem 3.14 follows along the lines of work
of Mazur~\cite{mazu:ripd} and Momose~\cite{momo:split} for elliptic curves
over $\Q$. However, at a key point we must introduce an analytic argument on
average special values of $L$-functions; it is because of the analytic
step that we can prove the main theorem of the paper only for $p \geq
211$.  It should be emphasized that the remaining cases can be handled
by a finite, though at present unfeasible, computation.  In particular,
the methods of this paper yield the following fact
(Proposition~\ref{pr:finitecomp})

\begin{intro}  Let $p>13$ be a prime.  Suppose there exists either
\begin{itemize}
\item a newform in $S_2(\Gamma_0(2p^2))$ with $w_p f = f$ and $w_2 f =
-f$; or 
\item a newform in $S_2(\Gamma_0(p^2))$ with $w_p f = f$,
\end{itemize}
such that $L(f \tensor \chi, 1) \neq 0$, where $\chi$ is the Dirichlet
character of conductor $4$.  Then the equation $A^4 + B^2
= C^p$ has no primitive non-trivial solutions.
\end{intro}

The author is grateful to Henri Darmon and Emmanuel Kowalski for
useful discussions about the theorems proved here, and to the referee
for thorough and helpful remarks.

\section{Twisted modular curves and their Jacobians}

A main theme of the present article is the analysis of certain twisted
versions of modular curves.  We begin with the ``untwisted''
versions.  Recall that, for any prime $p$, the curve $X^{split}(p)$ is
a coarse moduli space parametrizing pairs $(E,\set{A,B})$ where $E$ is
an elliptic curve and $\set{A,B}$ is an unordered pair of distinct
cyclic subgroups of $E[p]$.  Similarly, $X^{ns}(p)$ parametrizes pairs
$(E,N)$ where $N$ is a pair of Galois-conjugate points in $\P E[p]
\tensor_{\Fp} \F_{p^2}$.  We call such an $N$ a {\em non-split
structure} on $E$.  Note that $X^{split}(p)$ and $X^{ns}(p)$ are known
to have smooth models over $\Z[1/p]$.

\begin{defn} Let $p$ be a prime and $m$ a positive integer prime to
$p$.  We define
\beq
X_0^s(m;p) = X_0(m) \times_{X(1)} X^{split}(p)
\eeq
and
\beq
X_0^{ns}(m;p) = X_0(m) \times_{X(1)} X^{ns}(p).
\eeq
\end{defn}

The curves $X_0(mp), X_0^s(m;p)$ and $X_0^{ns}(m;p)$ have involutory
automorphisms $w_m$ arising from the Fricke involution $w_m$ on
$X_0(m)$. Let $K$ be a quadratic field.  For $X = X_0(mp),
X_0^s(m;p)$ or $X_0^{ns}(m;p)$, let $X^K/\Q$ be the twisted form of
$X$ admitting an isomorphism $\phi:X^K/K \ra X/K$ such that
$\phi^\sigma = w_m \circ \phi$. Note that
$X^K(\Q)$ can be described as the subset of $P \in X(K)$ such that
$P^\sigma = w_m P$, for $\sigma$ a generator of $\Gal(K/Q)$.

Recall that a {\em $\Q$-curve} $E/K$ is an elliptic curve which is
isogenous to its Galois conjugate.  If there exists such an isogeny of
degree $d$, we say $E/K$ is a $\Q$-curve of degree $d$.

In \cite[Prop.\ 2.3]{elle:qcurves} we define a
projective mod $p$ Galois representation
\beq
\P \bar{\rho}_{E,p}: \GalQ \ra \PGL_2(\F_p)
\eeq
associated to any $\Q$-curve.

We will later show that under certain circumstances, $\P
\bar{\rho}_{E,p}$ has large image.  Our main tool is the following
proposition. 


\begin{prop}  Let $d$ be a square-free positive integer, and let $E/K$
be a $\Q$-curve of degree $d$ over a quadratic number field $K$.  Let
$p$ be a prime not dividing $d$. 

Suppose the image of $\P \bar{\rho}_{E,p}$ lies in a
Borel subgroup of $\PGL_2(\F_p)$ (resp. the normalizer of a split
Cartan subgroup, normalizer of a non-split Cartan subgroup.)  Then the
point of $X_0(dp)(K)$ (resp. $X_0^{s}(d;p)(K)$, $X_0^{ns}(d;p)(K)$)
corresponding to $E$ is
a point of $X_0(dp)^K(\Q)$ (resp. $X_0^{s}(d;p)^K(\Q)$,
$X_0^{ns}(d;p)^K(\Q)$.
\label{pr:twistedx}
\end{prop}

\begin{proof}  We will discuss the case of $X_0(dp)$; the other two
are similar.  Let $C_p$ be the cyclic subgroup of $E[p]/K$ which is
fixed by $\GalK$.  Let $\mu: E^\sigma \ra E$ be the degree $d$
isogeny, and let $C_d = \mu(E^\sigma[d])$.  Now $\P
\rho_{E,p}(\sigma)$ acts on $\P E[p]$ by sending $x$ to
$\mu(x^\sigma)$.  Since, by hypothesis, $\P \rho_{E,p}(\sigma)$ fixes $C_p$, we have that $C_p^\sigma = \mu^{-1} C_p$.  So
indeed
\beq
(E,C_d,C_p)^\sigma = (E/C_d, E[d]/C_d, \mu^{-1} C_p) = w_d(E,C_d,C_p),
\eeq
which proves the desired result.
\end{proof}

\section{Good reduction for $\Q$-curves with rational level structures}

In this section we use arguments derived from Mazur's foundational
paper~\cite{mazu:ripd} in order to show that points on twisted modular
curves of large level must have good reduction at all large primes.

\begin{prop}  Let $R$ be a finite extension of $\Z_\ell$ with fraction
field $L$ and maximal ideal $\lambda$, and let $X/R$ be a stable
curve.  Write $X^{smooth}$ for the smooth part of $X$.  Suppose $A/R$ is a
semi-abelian scheme with a morphism $\phi: X^{smooth} \ra A$.  Let $x$ and
$y$ be distinct sections in $X^{smooth}(R)$, such that $\phi(x) = 0$.  Write
$x_\lambda$ and $y_\lambda$ for the restrictions of $x$ and $y$ to the
closed fiber of $\Spec R$.

 Then suppose
\begin{itemize}
\item The absolute ramification index $e_\ell$ of $R$ is less than $\ell-1$;
\item $\phi$ is a formal immersion at $x_\lambda$.
\item The restriction of $\phi(y)-\phi(x)$ to $A(L)$ is of finite order.
\end{itemize}
Then $x_\lambda$ and $y_\lambda$ are distinct points of
$X^{smooth}(R/\lambda)$.
\label{pr:formimm}
\end{prop}

\begin{proof}
The proposition is essentially Corollary $4.3$ of \cite{mazu:ripd}.
For the reader's convenience, we recount the argument here.  Let $q
\in A(R)$ be the point $\phi(y) - \phi(x)$.

First of all, suppose that the restriction to $q$ to $A(L)$ is a non-zero
torsion point of order $m$.  Then it follows from a theorem of
Raynaud~\cite[Prop.\ 1.1]{mazu:ripd} that the specialization
of $q$ to the closed point of $\Spec R$ also
has exact order $m$.  In particular, $q$ does not reduce to $0$
mod $\lambda$.  Thus, $y$ does not reduce to $x$ mod $\lambda$.

So we may assume that $\phi(y) = \phi(x)$.  Suppose $x_\lambda = y_\lambda$.
The section 
\beq
x: \Spec R \ra X
\eeq
restricts to a map from $\Spec R$ to the spectrum of the completed local ring
$\hat{\OO}_{X,x_\lambda}.$  So $x$ yields a map of rings
\beq
\tilde{x}: \hat{\OO}_{X,x_\lambda} \ra R. 
\eeq
Likewise, $y$ yields a map
\beq
\tilde{y}: \hat{\OO}_{X,y_\lambda} = \hat{\OO}_{X,x_\lambda} \ra R.
\eeq
Since $x \neq y$, the two maps $\tilde{x}$ and $\tilde{y}$ are distinct. 
Let $0_\lambda$ be the identity in $A(R/\lambda)$.  Then the fact that
$\phi$ is a formal immersion along $x$ means precisely that the map of
completed local rings
\beq
\tilde{\phi}: \hat{\OO}_{A,{0_\lambda}} \ra \hat{\OO}_{X,x_\lambda}
\eeq
is a surjection.  But this contradicts the fact that $\tilde{x} \circ
\tilde{\phi} = \tilde{y} \circ \tilde{\phi}.$
\end{proof}

\begin{prop} Let $K$ be a quadratic field, and $E/K$ be a $\Q$-curve
of squarefree degree $d$.  Suppose $\P \bar{\rho}_{E,p}$ is reducible for some
$p = 11$ or $p > 13$  with $(p,d) = 1$. Then $E$ has potentially good
reduction at all primes of $K$ of characteristic greater than $3$. 
\label{pr:borel}
\end{prop}

\begin{proof}  Let $C/\Q$ be the twisted modular curve $X_0(dp)^K$, and
take $J = \Jac(C)$.  By Proposition~\ref{pr:twistedx},
$E$ corresponds to a point $P \in C(\Q)$. We will think
of $P$ as lying in $C(K) = X_0(dp)(K)$.  Let $\lambda$ be a prime of
$K$ of characteristic greater than $3$.  We denote by $\infty$ the
usual cusp of $X_0(dp)$, and by $\infty_p$ the corresponding cusp of $X_0(p)$.

Define a map $g: X_0(dp)\ra J_0(p)$ by the rule
\beq
g((E,C_d,C_p)) = [E,C_p] + [E/C_d, C_p/C_d] - 2[\infty_p].
\eeq

Let $f$ be a Hecke eigenform in $S_2(\Gamma_0(p))$.  The map $\pi_f
\circ g: X_0(dp) \ra A_f$ extends by the N\'{e}ron mapping property to
a map, also denoted $\pi_f \circ g$, from the smooth part of
$X^0(dp)/\Z$ to the N\'{e}ron model of $A_f$, which we also denote $A_f$.

\begin{lem} The map $\pi_f \circ g: X_0^{smooth}(dp)/\Z[1/2] \ra A_f/\Z[1/2]$
is a formal immersion at the point $\infty \in
X_0^{smooth}(dp)(\F_\lambda)$.
\label{le:borelimm} 
\end{lem}

\begin{proof}
First of all, $\pi_f \circ g$ induces an isomorphism of the residue
fields attached to the points $\infty \in X_0^{smooth}(dp)(\F_\lambda)$
and $0 \in A_f(\F_\lambda)$, both of which are $\F_\lambda$.  So it suffices to
show that the map
\beq
(\pi_f \circ g)^*: \Cot_0(A_f/\F_\lambda) \ra \Cot_\infty(X_0)(dp)/\F_\lambda
\eeq
is a surjection~\cite[17.4.4]{grot:egaiv4}.

Let $g_1: X_0^{smooth}(dp)/\Z[1/2] \ra J_0(p)/\Z[1/2]$ be the morphism
sending the point $(E,C_d,C_p)$ to the divisor $[(E,C_p)] - [\infty]$.
Likewise, let $g_d$ be the morphism sending $(E,C_d,C_p)$ to
$[(E/C_d,C_p/C_d)] - [\infty]$.  Then $g = g_1 + g_d$.  Now $\pi_f \circ
g_1$ factors as
\beq
X_0(dp) \stackrel{p}{\ra} X_0(p) \stackrel{j}{\ra} A_f,
\eeq
where $p$, the ``forgetting of $d$-structure'' morphism, is unramified
at $\infty$ and $j$ is a formal immersion at the
point at infinity of $X_0(p)(\F_\lambda)$ (\cite[Prop\ 3.1]{mazu:ripd}.)  It
follows that the composition $\pi_f \circ g_1$ is a formal immersion
at $\infty$, whence $(\pi_f \circ g_1)^*$ is a surjection on
cotangent spaces.  On the other hand, $\pi_f \circ g_d$ factors as
\beq 
X_0(dp) \stackrel{p \circ w_d}{\ra} X_0(p) \stackrel{j}{\ra} A_f
\eeq
But $p \circ w_d$ is ramified at $\infty$, so the image of the
map $(p \circ w_d)^*$ in $\Cot_{\infty}(X_0)(dp)/\F_\lambda$ is zero.
Therefore, $(\pi_f \circ g_d)^*$ is also the zero map on cotangent
spaces.  So $(\pi_f \circ g)^* = (\pi_f \circ g_1)^*$ is surjective,
as desired.
\end{proof}

By the hypotheses on $p$, we can, and do, choose an eigenform $f \in
S_2(\Gamma_0(p))$ such that $A_f(\Q)$ is a finite group.  For instance,
we may choose $f$ so that $A_f$ lies in the Eisenstein quotient of
$J_0(p)$.  We now want to derive a contradiction from \ref{pr:formimm}. 

Suppose $E$ has potentially multiplicative reduction at $\lambda$.
Applying an Atkin-Lehner involution if necessary, we may assume that
$P$ reduces to $\infty$ mod $\lambda$.

We will now apply Proposition~\ref{pr:formimm}, where  $L = K_\lambda,
X = X_0(dp),$ and $A = A_f$.  Take $x$ to be the section
$\infty$, and $y$ to be the section given by $P$.  It follows from the
potentially multiplicative reduction of $E$ that $y$ is a section of
$X_0(dp)^{smooth}/\OO_{K_\lambda}$.

Finally, take $\phi$ to be the map $\pi_f \circ g.$  Since $L$ is a
quadratic extension of $\Q_\ell$, and $\ell > 3$, we have $e_\ell <
\ell - 1$.  By Lemma~\ref{le:borelimm}, $\phi$ is a formal immersion
at $x_\lambda = \infty_\lambda$.  

Let $\sigma$ be the nontrivial element of $\Gal(K/\Q)$;  then 
\beq
g(P)^\sigma  
= [p(P^\sigma)] + [p(w_d P^\sigma)] - 2[\infty_p] =
[p(w_d P)] + [p(P)] - 2[\infty_p] = g(P).
\eeq
In particular, $\phi(y) = \pi_f(g(P))$ lies in $A_f(\Q)$, and
therefore has finite order.  Note also that $\phi(x) = \pi_f(g(\infty)) = 0$.

We can conclude that $x_\lambda \neq
y_\lambda$, contradicting the hypothesis that $P$ reduces to $\infty$
mod $\lambda$.
\end{proof}

\begin{prop} Let $K$ be a quadratic field, and $E/K$ be a $\Q$-curve
of squarefree degree $d$.  Suppose the image of $\P \bar{\rho}_{E,p}$
lies in the normalizer of a split Cartan subgroup of $\PGL_2(\F_p)$, for
$p = 11$ or $p > 13$ with $(p,d) = 1$. Then $E$ has good reduction at
all primes of $K$ not dividing $6$.
\label{pr:split}
\end{prop}

\begin{proof}
This case is very similar to that of Proposition~\ref{pr:borel}.
Here, the key lemma on formal immersions is due to
Momose~\cite{momo:split}.

Let $C$ be the twisted modular curve $X_0^s(d;p)^K$.  Then $E$
yields a point 
\beq
P \in C(\Q) \subset C(K) = X_0^s(d;p)(K).
\eeq

To be more precise, $P \in C(K)$ parametrizes the triple
$(E,C^0_d,\set{A^0_p,B^0_p})$ where $C^0_d$ is the kernel of the
isogeny between $E$ and its Galois conjugate, and $\set{A^0_p,B^0_p}$
is the chosen pair of cyclic subgroups fixed by the action of $\GalK$.

Let $X_0^{sc}(d;p)$ be the modular curve parametrizing quadruples
$(E,C_d,A_p,B_p)$ where $C_d$ is a cyclic subgroup of $E$ of order $d$
and $A_p$ and $B_p$ are linearly independent cyclic subgroups of $E$
of order $p$.  Similarly, let $X^{sc}(p)$ be the curve parametrizing
triples $(E,A_p,B_p)$.  Note that $X_0^{sc}(d;p)$ (resp.\ $X^{sc}(p)$)
is naturally a double cover of $X_0^s(d;p)$ (resp.\ $X^{split}(p)$).
Write $w_p$ for the involution of $X_0^{sc}(d;p)$ switching $A_p$ and
$B_p$.  A cusp of $X_0^{sc}(d;p)$ is determined by its image in
$X_0(d)$ and $X^{sc}(p)$.  We write $(c,c')$ for the cusp of
$X_0^{sc}(d;p)$ whose image in $X_0(d)$ is $c$ and whose image in
$X^{sc}(p)$ is $c'$.

Let
\beq
h: X_0^{sc}(d;p)/{\Q} \ra J_0(p)/{\Q} 
\eeq
be the map defined by
\beq
h((E,C_d,A_p,B_p)) = [(E,A_p)] - [(E/B_p,E[p]/B_p)]
+ [(E/C_d,A_p/C_d)] - [(E/(B_p+C_d),E[p]/(B_p + C_d))].
\eeq
Then we have a diagram of schemes over $\Q$
\beq
\begin{CD}
X_0^{sc}(d;p) @>>> X_0^{s}(d;p)\\
@VhVV                @Vh^{-}VV   \\
J_0(p) @>>> (1-w_p)J_0(p) 
\end{CD}
\eeq
where $h^{-}$ is defined to make the above diagram commute.

The moduli problem over $\Spec \Q$ coarsely represented by
$X_0^{sc}(d;p)$ can be extended to a moduli problem over $\Spec \Z$ as
in \cite{katz:kama}; this moduli problem is coarsely represented by a
curve over $\Spec \Z$, which we also denote $X_0^{sc}(d;p)$.  (Note
that $X_0^{sc}(d;p)$ is isomorphic to $X_0(dp^2)$.)  The curve $X_0^{sc}(d;p)$
is smooth away from characteristics dividing $dp$.  In characteristics
dividing $d$, the reduction of $X_0^{sc}(d;p)$ is smooth away from
supersingular points, and in particular is smooth at all cusps.  In
characteristic $p$, the reduction of $X_0^{sc}(d;p)$ is made up of
three components, parametrizing triples $(E,C_d,A_p,B_p)$
where, respectively:
\begin{itemize}
\item $A_p \cong \mu_p$ \'{e}tale-locally;
\item $B_p \cong \mu_p$ \'{e}tale-locally;
\item neither $A_p$ nor $B_p$ is \'{e}tale-locally isomorphic to $\mu_p$.
\end{itemize}

The smooth part of the special fiber at $p$ contains the ordinary
locus in the first two components.  (See \cite[\S 1]{momo:split} and
\cite[(13.5.6)]{katz:kama} for facts used here about the special
fiber.)

By the N\'{e}ron mapping property, $h$ extends to a map from
$X_0^{sc;smooth}(d;p)/\Z$ to the N\'{e}ron model of $J_0(p)$.

Let $f$ be a Hecke eigenform in $S_2(\Gamma_0(p))$ such that $w_p f =
-f$.  Then we have a projection map from $J_0(p)/\Z$ to 
$A_f/\Z$. 

We base change this map to $\Z[\zeta_p][1/2]$ to
obtain a map
\beq
\pi_f: J_0(p)/{\Z[\zeta_p][1/2]} \ra A_f/{\Z[\zeta_p][1/2]}.
\eeq

\begin{lem} Let $\lambda$ be a prime of $\Q(\zeta_p)$ with residue field
  $\F_\lambda$.  The composition $\pi_f \circ h$ is a formal immersion
at every cusp $(\infty,c')$ of $X_0^{sc;smooth}(d;p)_{\F_\lambda}$, for all
$\lambda \notdiv 2p$.  For all $\lambda \notdiv 2$, the composition
$\pi_f \circ h$ is a formal immersion at the cusp $(\infty,\infty)$ of $X_0^{sc;smooth}(d;p)_{\F_\lambda}$.
\label{le:splitimm}
\end{lem}


\begin{proof}
We proceed much as we did in Lemma~\ref{le:borelimm}.  First, note
that the hypotheses on $\ell$ guarantee that $(c,c')$ lies in the
smooth locus of $X_0^{sc}(d;p)$ as in \cite[\S 1]{momo:split}.

Each cusp of $X_0^{sc}(d;p)_{\F_\lambda}$ is defined over $\F_\lambda$
(\cite[proof of (2.5)]{momo:split}), and $\pi_f \circ h$ is defined
over $\F_\lambda$ by definition; it follows that for each cusp $c$,
the residue field of $c$ and the residue field of $\pi_f
\circ h(c)$ are both $\F_\lambda$.  So in order to prove the lemma, it
suffices to show that the map
\beq
(\pi_f \circ h)^*: \Cot_0(A_f/\F_\lambda) \ra \Cot_{(\infty,c')}X_0^{sc}(d;p)/\F_\lambda
\eeq
is a surjection.

We now write $h = h_1 + h_d$,
where
\beq
h_1((E,C_d,A_p,B_p)) = [(E,A_p)] - [(E/B_p,E[p]/B_p)]
\eeq
and $h_d = h_1 \circ w_d$.  Now $\pi_f \circ h_1$ factors as
\beq
X_0^{sc;smooth}(d;p) \stackrel{a}{\ra} X^{sc;smooth}(p) \stackrel{\pi_f \circ
g}{\ra} A_f.
\eeq 
Here $g$ is the morphism from $X^{sc;smooth}(p)$ to $J_0(p)$ defined by
Momose~\cite[\S 2]{momo:split}, and $a$ is the ``forgetting of
$d$-structure'' morphism.  Again, $a$ is unramified at $(\infty,c')$, and
$\pi_f \circ g$ is a formal immersion at the cusp $c'$ of
$X^{sc;smooth}(p)$, by \cite[(2.5)]{momo:split}.  So $(\pi_f \circ
h_1)^*$ is a surjection on cotangent spaces. 
 
On the other hand, $w_d \circ a$ is ramified at
$(\infty,c')$, so $(\pi_f \circ h_d)^*$ is the zero map on cotangent
spaces.  So $(\pi_f \circ h)^*$ is surjective, as desired.
\end{proof}

We now make the further stipulation on $f$ that $A_f(\Q)$ is a finite
group.  (Again, we may choose $f$ 
such that $A_f$ is a quotient of the Eisenstein quotient of $J_0(p)$.) 
Let $\sigma$ be the nontrivial element of $\Gal(K/\Q)$.  We have
\beq
(E,C^0_d,\set{A^0_p,B^0_p})^\sigma = (E/C^0_d, E[d]/C^0_d,
\set{A^0_p/C^0_d,B^0_p/C^0_d}). 
\eeq
It follows immediately that $h^{-}(P^\sigma) = h^{-}(P)$.  Let $Q \in
X_0^{sc}(d;p)(\bar{\Q})$ be a point lying over $P$.  Then $h(Q) =
h^{-}(P) \in [(1-w_p)J_0(p)](\Q)$, so $\pi_f(h(Q))$
lies in $A_f(\Q)$, and is thus of finite order.

Let $M$ be the field of definition of $Q$.  Let $G$ be the group of
automorphisms of $X_0^{sc}(d;p)$ generated by $w_d$ and $w_p$; then $G
\cong (\Z/2\Z)^2$ and $Q$ lies over a $\Q$-point of $X_0^{sc}(d;p)/G$,
whence $M$ is a subfield of a biquadratic field over $\Q$.  Suppose
$\lambda$ is a prime of $M$ such that $E$ has potentially
multiplicative reduction at $\lambda$.  Write $(c,c') \in
X_0^{sc}(d;p)(\bar{\Q})$ for the cusp to which $Q$ reduces mod
$\lambda$.  Applying $w_d$ if necessary, we may assume that
$c=\infty$.

If $\lambda \notdiv p$, the map $\pi_f \circ h$ is a formal immersion
at $(c,c')_\lambda$ by Lemma~\ref{le:splitimm}.  Suppose $\lambda |
p$.  Since $E$ has potentially multiplicative reduction, it acquires
semistable reduction after a quadratic extension $M'$ of $M$; since
$M$ is biquadratic, the absolute ramification index of $M'$ at $p$ is
at most $4$.  If $Q$ reduces to a cusp other than $0$ and $\infty$,
the group schemes $A^0_p$ and $B^0_p$ are both \'{e}tale over
$\OO_{M'}$ (\cite[proof of (2.5)]{momo:split} whence, by Weil pairing,
$\mu_p$ is \'{e}tale over $\OO_{M'}$; this makes the absolute
ramification index of $M'$ over $p$ at least $p-1$, a contradiction.
If $Q$ reduces to $0$, we can act on $Q$ by $w_p$ to make $c' =
\infty$.  Now, by Lemma~\ref{le:splitimm}, the map $\pi_f \circ h$ is
a formal immersion at $(c,c')_\lambda$.

We will now apply Proposition~\ref{pr:formimm}, using $L =
M_{\lambda}, X = X_0^{sc}(d;p),$ and $A = A_f$.  Take $x$ to be the
cuspidal section $(c,c')$, and $y$ to be the section given by $Q$.
Finally, take $\phi$ to be the map $\pi_f \circ h$.  Since $M$ is a
subfield of a biquadratic extension of $\Q$, its ramification degree
is at most $2$ over any odd prime.  So we have $e_\ell <
\ell - 1$.  Now the conclusion of
Proposition~\ref{pr:formimm} contradicts the hypothesis that $Q$ and
$(c,c')$ reduce to the same point of $X_0^{sc}(d;p)(M/\lambda)$.
\end{proof}

\medskip

We now turn to the case of $\Q$-curves $E$ whose mod $p$ Galois
representations have image in the normalizer of a non-split Cartan
subgroup.  This case is more difficult, due to the absence of rank $0$
quotients of $J^{ns}(p)$.  However, we show below by analytic means
that the Jacobian of the twisted modular curve $X_0^{ns}(d;p)^K$ does
have rank $0$ quotients; we then obtain a good reduction theorem on
$E$ using a formal immersion result of Darmon and
Merel~\cite{darm:wind}.


\begin{prop} Let $K$ be an imaginary quadratic field, and $E/K$ be a $\Q$-curve
of squarefree degree $d$.  There exists a constant $M_{d,K}$ with the following property.

Suppose the image of $\P \bar{\rho}_{E,p}$ lies in the
normalizer of a non-split Cartan subgroup of $\PGL_2(\F_p)$, for
$p > M_{d,K}$. Then $E$
has potentially good reduction at all primes of $K$.
\label{pr:ns}
\end{prop}

\begin{rem}  All is not lost if $K$ is a real quadratic field.  We
will see below that when $K$ is imaginary, there exists a newform $f$
on level $p^2$ satisfying the conditions of
Proposition~\ref{pr:analytic}.  When $K$ is real, there still may be
newforms $f$ on other levels satisfying those conditions.  However,
the methods described here cannot treat the case of arbitrary $K$ and
$d$.  For instance, suppose $d$ is prime, and $K$ is a real quadratic field
in which $d$ is inert.  Then if $f$ is either a newform in
$S_2(\Gamma_0(p^2))$ with $w_p f = f$, or a newform in
$S_2(\Gamma_0(dp^2))$ with $w_p f = f$ and $w_d f = -f$, we have by a
theorem of Weil \cite[Th.\ 6]{li:newforms} that $f \tensor \chi_K$ has
negative functional equation, so that $L(f \tensor \chi_K, 1) = 0$.
This leaves us in much the same position as one who
tries to use Mazur's method to control points on $X^{ns}(p)$; the
Jacobian of the curve in question, assuming Birch-Swinnerton-Dyer, has
no rank $0$ quotient.
\end{rem}
 
\begin{proof}

Let $C$ be the twisted modular curve $X_0^{ns}(d;p)^K$, and write $J =
\Jac(C)$.  Then $E$ yields a point $P \in C(\Q)$.

We have an isomorphism $J \times_\Q K \cong J_0^{ns}(d;p)/K$; from a
result of Chen and Edixhoven~\cite{chen:jacobians},\cite{desm:chen}
there is a surjective homomorphism
\beq
\alpha: J_0^{ns}(d;p) \ra J'_0(dp^2)/w_p.
\eeq
where $J'_0(dp^2)$ is the $p$-new quotient of $J_0(dp^2)$.  
  
Let ${\mathbf T}$ be the algebra generated by Hecke operators of degrees
prime to $dp$, together with the group $W$ of Atkin-Lehner involutions
of degrees dividing $d$.  It follows from Theorem 2 of
\cite{desm:chen} that the map $\alpha$ is compatible with the
action of ${\mathbf T}$ on either side (in that theorem, let $\CC$ be the
isogeny category of abelian varieties endowed with an action of
${\mathbf T}$, and $M$ the Jacobian of $X_0(d) \times_{X(1)} X(p)$.)

Suppose $f$ is either
\begin{itemize}
\item a newform in $S_2(\Gamma_0(dp^2))$ with $w_p f = f$ and $w_d f = -f$;
\item a newform in $S_2(\Gamma_0(d'p^2))$ with $d'$ a proper divisor
of $d$ and $w_p f = f$.
\end{itemize}
In each case, we have a quotient morphism
\beq
\pi_f: J_0(dp^2) \ra A_f
\eeq
such that the action of $w_d$ on $J_0(dp^2)$ induces the involution
$-1$ on $A_f$.  In case $d$ is not prime, we can and do choose our
$\pi_f$ such that the quotient $A_f$ is preserved by the whole group $W$
of Atkin-Lehner involutions.  More precisely:  for each $e | (d/d')$,
we have a map $B_e: J_0(d'p^2) \ra J_0(dp^2)$.  Choose a character $\chi:
(\Z/(d/d')\Z)^* \ra \pm 1$.  It follows from Lemma 26 of
\cite{atki:al} that the quotient map
\beq
I_\chi: \sum_{e | (d/d')} e B_e \chi(e) : J_0(d'p^2) \ra J_0(dp^2)
\eeq
has image stable under the action of $W$, and on which the action of
$w_d$ on the quotient is $w_{d'}$ twisted by the scalar $\chi(d)$.  So
if we choose $\chi$ such that $\chi(d)$ and the eigenvalue of $w_{d'}$ on $A_f$
have opposite signs, then the image under $I_\chi$ of $A_f \subset
J_0(d'p^2)$ is a subvariety of $J_0(dp^2)$, isogenous to $A_f$, which
is stable under $W$ and on which $w_{d}$ acts as $-1$.  Now let
$\pi_f$ be projection onto that subvariety.

Composing with $\alpha$ yields a morphism from $J_0^{ns}(d;p)$ to
$A_f$; replacing $A_f$ by an isogenous variety $A'_f$, we have a quotient 
\beq
\pi'_f: J_0^{ns}(d;p) \ra A'_f
\eeq
which is compatible with the action of ${\mathbf T}$ and has connected
kernel.  In particular, $\pi'_f \circ w_d  = -\pi'_f$.

Denote $\chi_K$ by $\chi$, and write $A_f \tensor \chi$ for the twist
of $A_f$ by $\chi$. Let $\sigma$ be the non-trivial element of $\Gal(K/\Q)$.
Then we have a commutative diagram of abelian varieties
over $K$:
\beq
\begin{CD}
J_0^{ns}(d;p) @>{\pi'_f}>> A'_f \\
@ViVV @VjVV \\
J @>>> A'_f \tensor \chi
\end{CD}
\eeq

where $i$ and $j$ are isomorphisms such that $i^\sigma =  w_d \circ i$ and
$j^\sigma = -j$.  Let
\beq
\psi_f: J \ra A'_f \tensor \chi
\eeq
be the composition $j \circ \pi'_f \circ i^{-1}$.  Then
\beq
\psi_f^\sigma = j^\sigma \circ (\pi'_f)^\sigma \circ (i^{-1})^\sigma
= (-j) \circ \pi'_f \circ w_d \circ i^{-1}
= (-j) \circ (-\pi'_f) \circ i^{-1} = \psi_f.
\eeq
In other words, $\psi_f$ is defined over $\Q$.

Let $R_0$ be the ring of integers of the number field $K(\zeta_p +
\zeta_p^{-1})$, and let $R = R_0[1/2dp]$.  Then $X_0^{ns}(d;p)$ has a
smooth model over $R$ and the cusp $\infty$ of
$X_0^{ns}(d;p)$ is defined over $R$~\cite[\S 5]{darm:wind}. 


We can define a map
\beq
h: X_0^{ns}(d;p)/R \ra J_0^{ns}(d;p)/R
\eeq
by setting $h(P) = [P] - [\infty]$.

\begin{lem} Let $\lambda$ be a prime of $R$.  Then
the map
\beq
\pi'_f \circ h: X_0^{ns}(d;p)/R \ra A'_f/R
\eeq
is a formal immersion at the point $\bar{\infty}$ of 
$X_0^{ns}(d;p)(\F_\lambda)$.
\label{le:nsformimm}
\end{lem}

\begin{proof}
This fact is almost precisely Lemma 8.2 of \cite{darm:wind}.

One difference is that our quotient $A'_f$ is not preserved by $T_n$
for all $n$ prime to $p$, but by $T_n$ for all $n$ prime to $dp$ and all $w_{d'}$ for $d'
| d$.  We need to prove that, as in \cite{darm:wind}, there exists a
differential form $\omega$ on $A'_f$ whose associated modular form $\tilde{g}
= \sum a_n(\tilde{g}) q^{n/p}$ has $a_1(\tilde{g}) \neq 0 \pmod
\lambda$.  In fact, we will prove this for any quotient $A$ of
$J_0^{ns}(d;p)$ which is preserved by ${\mathbf T}$ and which is not
killed by $\alpha$.

Write $S(A)$ for the vector space of weight $2$ cusp forms attached to
$A$.  Suppose that $a_1(\tilde{g}) = 0 \pmod \lambda$ for
every $\tilde{g}$ in $S(A)$.

Choose some $\tilde{g}$ in $S(A)$ which is an eigenform for ${\mathbf
T}$, which is not in the kernel of $\alpha$, and which does not reduce
to $0$ mod $\lambda$.  The form 
\beq
g = \sum a_n(\tilde{g}) q^n
\eeq
is a form on $\Gamma_1(dp^2)$ which is also an eigenform for
${\mathbf T}$.  (We remark, however, that the Hecke eigenform on
$\Gamma_0(dp^2)$ associated to $\tilde{g}$ via the map $\alpha$ does
not necessarily have the same Hecke eigenvalues as $g$.)

Let $g_0$ be a newform on some level $dp^2/M$ with the
same eigenvalues as $g$.  Then we can write
\beq
g = \sum_{d'|M} \alpha_{d'} B_{d'} g_0.
\eeq
where $B_d$ is the Hecke operator sending $f(\tau)$ to $f(d\tau)$.
(In this and all other discussion of Hecke operators, we follow the
notation of \cite{atki:al}.)

Suppose $\alpha_{e} \neq 0 \pmod \lambda$ for some $e|M$ with $(e,p)
= 1$.  By \cite[Prop.\ 1.5]{atki:twistw},
\beq
w_{e} g = c \alpha_{e} B_1 g_0 + \sum_{d'|M, d' > 1} B_{d'} h_{d'}
\eeq
where $c$ is a constant not divisible by $\lambda$ and the $h_{d'}$
are other modular forms.  Since $a_1(B_{d'} h_{d'}) = 0$ for all
$d'>1$, we see that $a_1(w_e g) \neq 0 \pmod \lambda$, a
contradiction.  We conclude that $\alpha_{d'} = 0 \pmod \lambda$
unless $p|d'$; so $g$ is congruent $\pmod \lambda$ to a form in the
image of $B_p$, which implies that $a_n(g) = 0 \pmod \lambda$ unless
$p|n$.  This in turn implies that $\tilde{g}$ is fixed by the action
not only of the normalizer of nonsplit Cartan in $\GL_2(\Z/p\Z)$, but
of a Borel subgroup as well.  So $\tilde{g}$ is a modular form on
$\Gamma_0(d)$ and is therefore killed by $\alpha$, a contradiction.
\end{proof}

If $f$ has the Fourier expansion $\sum a_n q^n$, write $f \tensor \chi$
for the modular form $\sum \chi(n) a_n q^n$.  Then $f \tensor \chi$ is
a newform of some level $N$, and in particular there is an associated quotient
$A_{f \tensor \chi}$ of $J_0(N)$.  Moreover, the abelian varieties $A'_f
\tensor \chi$ and $A_{f \tensor \chi}$ are isogenous over $\Q$.

\begin{prop}  Suppose $K$ is an imaginary quadratic field, and $\chi =
\chi_K$ is the associated Dirichlet character.
For all sufficiently large $p$, there exists a weight $2$ cusp form
$f$, which is either
\begin{itemize}
\item a newform in $S_2(\Gamma_0(dp^2))$ with $w_p f = f$ and $w_d f = -f$;
\item a newform in $S_2(\Gamma_0(d'p^2))$ with $d'$ a proper divisor
of $d$ and $w_p f = f$,
\end{itemize}
such that $A_{f \tensor \chi}(\Q)$ is a finite group. 
\label{pr:analytic}
\end{prop}

We first explain how to finish the proof of Proposition~\ref{pr:ns}
assuming the result of Proposition~\ref{pr:analytic}. 

Suppose $E/K$ is a $\Q$-curve of degree $d$, meeting the
hypotheses of Proposition~\ref{pr:ns}.

First of all, suppose
$\lambda$ is a prime of $K$ dividing $p$.  If the reduction of $E$ at
$\lambda$ is potentially multiplicative, then the image of the
decomposition group $G_\lambda$ under $\P\bar{\rho}_{E,p}$ lies in a
Borel subgroup.  On the other hand, by hypothesis this image lies in
the normalizer of a non-split Cartan subgroup.  We conclude that the
size of this image has order at most $2$, which means that $K_\lambda$
contains $\Q(\zeta_p + \zeta_p^{-1})$.  This is impossible once $p
\geq 7$.

Now suppose that $E$
has potentially multiplicative reduction over a prime $\ell$ not
dividing $p$. 

The cusps of $X_0^{ns}(d;p)$ have minimal field
of definition $\Q(\zeta_p + \zeta_p^{-1})$~\cite[\S 5]{darm:wind},
and $K$ is linearly disjoint from $\Q(\zeta_p + \zeta_p^{-1})$; it
follows that the cusps of $X_0^{ns}(d;p)$ which lie over $\infty \in
X_0(d)$ form a single orbit under
$\GalK$. If $\lambda$ is a prime of $K(\zeta_p + \zeta_p^{-1})$ over
$\ell$, then the point $P \in X_0^{ns}(d;p)(K)$ parametrizing $E$
reduces mod $\lambda$ to some cusp $c$.  By applying Atkin-Lehner
involutions at the primes dividing $d$, we can ensure that $P$ reduces
to a cusp which lies over $\infty$ in $X_0(d)$.  By the transitivity
of the Galois action, we can choose $\lambda$ so that $P$ actually
reduces to the cusp $\infty$ mod $\lambda$.  Note that in order for a
$K$-point of $X_0^{ns}(d;p)$ to reduce to $\infty$, the residue field
$\OO_K/\lambda$ must contain $\zeta_p + \zeta_p^{-1}$; this implies
that $\ell^4 \equiv 1$ mod $p$, and in particular $\ell
\neq 2,3$ when $p \geq 7$.  Moreover, if $p$ is large enough, we have
$(d,\ell) = 1$.

Now take $f$ to be a form satisfying the conditions specified in
Proposition~\ref{pr:analytic}.  We have defined above a map 
\beq
X_0^{ns}(d;p)^{K}/R \ra A'_f \tensor \chi,
\eeq 
which is a formal immersion at $\bar{\infty}$ by
Lemma~\ref{le:nsformimm}.  Apply Proposition~\ref{pr:formimm} with
$X = X_0^{ns}(d;p)/R_\lambda, x = \infty, y = P,$ and $\phi = \pi'_f
\circ h$.  We first apply the argument of \cite[Lemma 8.3]{darm:wind} to
show that the point $\phi(P)$ is torsion in $A'_f(L)$,
where $L = K(\zeta_p + \zeta_p^{-1})$.  Let $n$ be an integer which
kills the subgroup of $J_0^{ns}(d;p)$ generated by cusps; such an $n$
exists by the Drinfel'd-Manin theorem.  Let $\sigma$ be a generator for
$\Gal(L/K)$; then $(\sigma-1)h(P)$ is killed by $n$, so $nh(P)$ lies in
$J_0^{ns}(d;p)(K)$.  Let $\tau$ be an element of $\Gal(L/\Q)$ not
lying in $\Gal(L/K)$.  Then $P^\tau = w_d P$, and 
\beq
n\phi(P)^\tau = n[P^\tau] - n[\infty^\tau] = n[w_d P] - n[w_d \infty] =
n w_d \phi(P) = -n \phi(P).
\eeq
So $n \phi(P)$ lies in the subgroup of $A'_f(K)$ on which $\tau$ acts
as $-1$.  By the hypothesis on $f$, this subgroup is finite.  We
conclude that $\phi(P)$ is torsion.

Since $\ell > 3$, the absolute ramification index of $R_\lambda$ at
$\ell$ is at most $2$. It now follows from
Proposition~\ref{pr:formimm} that $y$ and $x$ reduce to distinct
points of $X$, contradicting our hypothesis on $E$.

\medskip

It now remains only to prove Proposition~\ref{pr:analytic}.

We say a form $f$ is {\em $p$-new} if it is not in the space of old
forms arising from $S_2(\Gamma_0(dp))$.

By a theorem of Kolyvagin and Logachev~\cite{koly:kolo}, building on
results of Gross-Zagier, Bump-Friedberg-Hofstein, and Murty-Murty, it
suffices to show that there exists a weight $2$ newform $f$ on level
$p^2$ such that $w_p f = f$ and $L(f \tensor \chi,1) \neq 0.$

Our method will be to show the stronger statement that the (suitably
weighted) {\em average} value of the above $L$-function over a certain
class of forms is nonzero.

Let $\FF$ be a Petersson-orthogonal basis for
$S_2(\Gamma_0(p^2))$ such that each $f \in \FF$ is an
eigenform for all Hecke operators $T_\ell$ where $\ell$ is prime to
$p$, and for the Atkin-Lehner involution $w_p$.

We define an average
\beq
V(p) = \sum_{\stackrel{f \in \FF}{w_p f = f}} a_1(f) L(f \tensor
\chi, 1).
\eeq

First of all, note that if $f$ is a form with $w_p f = -f$, then the
functional equation of $L(f \tensor \chi,s)$ has sign $\chi(-1) = -1$
by \cite[Th.\ 6]{li:newforms}.  So in this case $L(f \tensor \chi,
1) = 0$.  It follows that
\beq
V(p) = \sum_{f \in \FF} a_1(f) L(f \tensor \chi, 1).
\eeq
 


We may think of a Fourier coefficient $a_n$ as a linear functional in
$\Hom(S_2(\Gamma_0(p^2)),\C)$.  Likewise, write $L_\chi$ for the
functional sending $f$ to $L(f \tensor \chi, 1)$.  Now the Petersson
inner product on $S_2(\Gamma_0(p^2))$ defines an inner product on the
dual space $\Hom(S_2(\Gamma_0(p^2)),\C)$, and the average we are
studying is
\beq
V(p) = (a_1, L_\chi)
\eeq
with respect to this dual inner product.

In general, if $\ell_1,\ell_2$ are linear functionals on the space of
cuspforms for some $\Gamma_0(N)$, we write $(\ell_1,\ell_2)_N$ for the
corresponding Petersson product.  If $V$ is a subspace of
$S_2(\Gamma_0(N))$, we write $(\ell_1,\ell_2)_V$ for the Petersson
product restricted to $V$.  Finally, if $M|N$, we write
$(\ell_1,\ell_2)_N^M$ for the contribution to $(\ell_1,\ell_2)_N$ of
those forms which are new on level $M$.

The value of $V(p)$ can be estimated using the Petersson formula, as
in \cite{duke:average}.  In particular, we will show that, for $p$
large enough, $V(p)$ is nonzero, and thus that $L(f \tensor \chi,1)$
is nonzero for some $f$ in $\FF$.

\begin{lem} $V(p) = 4\pi + O(p^{-2+\epsilon}).$
\label{le:vp}
\end{lem}

\begin{proof} Immediate from the main theorem of \cite{elle:duke}; to
be precise, we have in general that
\begin{equation}
(a_m, L_\chi)_N = 4\pi\chi(m) + O(N^{-1+\epsilon})
\label{eq:amlchi}
\end{equation}
with constants depending only on $m,\chi,$ and $\epsilon$.
\end{proof}

But Lemma~\ref{le:vp} is not enough for us, since we require that
there be a {\em $p$-new} form $f$ with $L(f \tensor \chi,1) \neq 0$.
We must therefore show that the contribution of the $p$-old forms to
$V(p)$ is close to $0$.  We give below a general bound for the
contribution of $p$-old forms to $(a_m,L_\chi)_{p^2}$. This will
require an argument somewhat more intricate, but no deeper, than the
Petersson estimate for $V(p)$.  The problem of bounding the
contribution of oldforms is treated in \cite{iwan:lowlying}, but only
in case the level is square-free.  We have recently learned that a
paper of Akbary~\cite{akba:oldforms} also bounds the contribution of
oldforms in a similar situation.

\begin{rem} We expect that, for arbitrary fixed $N$, the contribution
of $p$-old forms to $(a_m,L_\chi)_{Np}$ is $o((a_m,L_\chi)_{Np})$;
proving this, when a high power of $p$ divides $N$, seems rather
complicated.
\end{rem}

The space of $p$-old forms on $\Gamma_0(p^2)$ is orthogonal to
the space of $p$-new forms.  So we can decompose the inner product
$(a_m, L_\chi)$ as 
\beq
(a_m, L_\chi) = (a_m,L_\chi)^{p-new} + (a_m,L_\chi)^p_{p^2}
\eeq


We will show that $(a_m,L_\chi)^p_{p^2}$ approaches $0$ as $p$ grows.

\begin{lem}  Let $f$ be a weight $2$ newform on $\Gamma_0(p)$, let
$\lambda_p(f)$ be the eigenvalue of $W_p$ on $f$,
and let $V_f$ be the space of forms on $\Gamma_0(p^2)$ arising from
$f$.  Let $m$ be a positive integer prime to $p$.  Then
\beq
(a_m, L_\chi)_{V_f} = \frac{p}{p^2-1}[1 + p^{-1}\chi(p)\lambda_p(f)](f,f)^{-1}a_m(f)L_\chi(f).
\eeq
\label{le:fromdp}
\end{lem}

\begin{proof}
The space $V_f$ is spanned by $B_1 f$ and $p B_p f$.

Then
\begin{equation}
(a_m,L_\chi)_{V_f} = 
\left[\begin{array}{cc} a_m(B_1 f) & a_m(p B_p f) \end{array} \right]
A^{-1} \left[\begin{array}{c} L_\chi(B_1 f) \\
p L_\chi(p B_p f) \end{array} \right]
\label{eq:amlchivf2}
\end{equation}
where $A$ is the symmetric matrix defined by $A_{ij} = (p^i B_i f, p^j
B_j f)$.

It follows from the definition of Petersson product that
\beq
A_{11} = A_{22} = [\Gamma_0(p^2):\Gamma_0(p)](f,f) = 
p(f,f).
\eeq   
We will now show that
\begin{equation}
(B_1 f, p B_p f) = -\lambda_p(f)(f,f).
\label{eq:b1bp}
\end{equation}

Recall that we can write the Petersson product of two forms $f$ and
$g$ as
\beq
C \Res_{s=2} L(s,f \times g)
\eeq
where $C$ is a constant independent of $f$ and $g$, and $L(s,f \times
g)$ is the Rankin-Selberg $L$-function defined by analytic
continuation of the series 
\beq
L(s,f \times g) = \sum_{n=1}^\infty a_n(f) a_n(g) n^{-s}.
\eeq
(see \cite[\S 1.6]{bump:aut}.)

Now 
\begin{eqnarray*}
C^{-1}(B_1 f, p B_p f) & = & p \Res_{s=2} \sum_{n=1}^\infty a_n(f) a_n(B_p
f) n^{-s} \\
& = &
p \Res_{s=2} \sum_{n=1}^\infty a_{pn}(f)a_n(f)n^{-s}p^{-s} \\
& = &
- p^{-1} \lambda_p(f)\Res_{s=2} \sum_{n=1}^\infty [a_n(f)]^2 n^{-s} \\
& = & - p^{-1} \lambda_p(f) C^{-1} (B_1 f, B_1 f).
\end{eqnarray*}
We now have
\beq
A = (f,f)\mat{p}{-\lambda_p(f)}{-\lambda_p(f)}{p}.
\eeq
Note that $\lambda_p(f) = \pm 1$.  So 
\beq
A^{-1} = (f,f)^{-1} (p^2 - 1)^{-1} \mat{p}{\lambda_p(f)}{\lambda_p(f)}{p}.
\eeq
Note that $(B_pf) \tensor \chi = \chi(p)B_p(f \tensor \chi)$.
Moreover, if $g$ is any modular form,
\beq
L(B_p g, 1) = \int^\infty_0 B_pg(iy) dy = \int^\infty_0 g(ipy) dy =
(1/p)L(g,1).
\eeq
We conclude that $L_{\chi}(p B_p f) = \chi(p) L_\chi(f)$.  So 
\begin{eqnarray*}
(a_m,L_\chi)_{V_f} & = & p(p^2-1)^{-1}(f,f)^{-1}\left[\begin{array}{cc} a_m(B_1 f) &
0 \end{array} \right] \mat{1}{p^{-1}\lambda_p(f)}{p^{-1}\lambda_p(f)}{1}
\left[\begin{array}{c} L_\chi(f) \\
\chi(p) L_\chi(f) \end{array} \right] \\
& = &
p(p^2-1)^{-1}[1 + p^{-1}\chi(p)\lambda_p(f)](f,f)^{-1}a_m(f)L_\chi(f).
\end{eqnarray*}
\end{proof}

Now $(a_m,L_\chi)^{p}_{p^2}$ is the sum over newforms $f$ of level $p$ of
$(a_m,L_\chi)_{V_f}$, which by Lemma~\ref{le:fromdp} is equal to 
\beq
\frac{p}{p^2-1}[1 +
  p^{-1}\chi(p)\lambda_p(f)](f,f)^{-1}a_m(f)L_\chi(f)
= \frac{p}{p^2-1}[a_m(f) - p^{-1}\chi(p)a_{mp}(f)]L_\chi(f)(f,f)^{-1}.
\eeq
Summing the above quantity over a Petersson-orthogonal basis of
newforms for $\Gamma_0(p)$ yields
\beq
(a_m,L_\chi)^{p}_{p^2} =  \frac{p}{p^2-1}(a_m -
p^{-1}\chi(p)a_{mp},L_\chi)_p.
\eeq
Now $(a_m,L_\chi)_p$ is bounded as $p$ grows by
\eqref{eq:amlchi}.   So it suffices to show that
$(a_{mp},L_\chi)_p$ is $o(p^2)$.  We will prove a version of this fact
with explicit constants, since these will be needed in the sequel.

\begin{lem} Let $p$ be a prime, $m$ a positive integer, $\chi$ a
quadratic character of conductor $q$ prime to $p$.  Then
\beq
(a_{mp},L_\chi)_p \leq
2\sqrt{3}m^{1/2}d(m)(1-e^{-2\pi/q\sqrt{p}})^{-1}
(4\pi+16\zeta^2(3/2)\pi^2p^{-3/2}).
\eeq
\label{le:convexity}
\end{lem}

\begin{proof}
Let $\FF_{p}$ be a Petersson-orthogonal basis of weight $2$ cuspforms on
$\Gamma_0(p)$.  Suppose furthermore that each $f \in \FF_{p}$ is
a Hecke eigenform.

First of all, note that $|a_p(f)| = |a_1(f)|$, and $|a_{mp}(f)| \leq
m^{1/2}d(m)|a_1(f)|$ by Weil bounds.  Now $f \tensor \chi$ is a
newform on level $M = pq^2$. The
functional equation for $L(f \tensor \chi, s)$ tells us that, for any
positive real $x$,
\beq
L(f \tensor \chi, 1) = \sum_{n>0} \chi(n) |a_n(f)| n^{-1} e^{-2\pi n/ x}
+ \sum_{n>0} \chi(n) |a_n(w_M (f \tensor \chi))| n^{-1} e^{-2\pi n x / M}.
\eeq
Since $f \tensor \chi$ is a newform, we have $|a_n(w_M(f
\tensor \chi))| = |a_n(f)|$.  We now set $x = \sqrt{M}$ and obtain the
bound
\beq
|L(f \tensor \chi,1)| \leq 2 \sum_{n > 0} a_n(f) n^{-1} e^{-2 \pi n /
\sqrt{M}} \leq 2(\sum_{n > 0} n^{-1/2} d(n) e^{-2 \pi n / \sqrt{M}})
|a_1(f)|.
\eeq

The sum over $n$ is of length approximately $\sqrt{M}$, and so has
value of order at most $M^{1/4 + \epsilon}$.  Working the constants
out is slightly intricate, so we satisfy ourselves with a much cruder
bound.  Since $d(n) \leq \sqrt{3n}$, we have 
\beq
|L(f \tensor \chi,1)| \leq 2\sqrt{3}(1-e^{-2 \pi/\sqrt{M}})^{-1}|a_1(f)|,
\eeq
so
\begin{eqnarray*}
|(a_m, L_\chi)| & \leq & \sum_{f \in \FF} |a_m(f)||L(f \tensor \chi,
 1)| \\
& \leq & \sum_{f \in \FF} m^{1/2}d(m) 2\sqrt{3}(1-e^{-2
\pi/\sqrt{M}})^{-1} a_1(f)^2 \\
& = & 2\sqrt{3}m^{1/2}d(m)(1-e^{-2\pi/q\sqrt{p}})^{-1} (a_1,a_1)_p.
\end{eqnarray*}
Now by Lemma 4 of \cite{elle:duke} we have
\beq
|(a_1,a_1)_p| \leq 4\pi + 16\zeta^2(3/2)\pi^2p^{-3/2}.
\eeq
This yields the desired result.
\end{proof}

We have now proved that
$(a_m,L_\chi)^p_{p^2}$ approaches $0$ as $p$ goes to $\infty$. Therefore,
$(a_1, L_\chi)^{p-new}$ approaches $V(p)$ as $p$ grows; in particular,
\beq
(a_1,L_\chi)^{p-new} \neq 0
\eeq
for $p$ sufficiently large.  We have now proved 
Proposition~\ref{pr:analytic}, and therefore also
Proposition~\ref{pr:ns}.
\end{proof}

\bigskip

Suppose that $K$ is a quadratic field, and $E/K$ is a $\Q$-curve of
degree $d$. If $\P\bar{\rho}_{E,p}$ does not
surject onto $\PGL_2(\F_p)$, then the image of $\P\bar{\rho}_{E,p}$
is contained in a maximal subgroup of $\PGL_2(\F_p)$; that is to say,
the image is contained in either a Borel subgroup, the normalizer of a
Cartan subgroup, or an exceptional subgroup isomorphic to $A_4, S_4$,
or $A_5$.  For any given $K$, there are only finitely many $p$ for
which it is possible that $\P\bar{\rho}_{E,p}$ has image contained in
an exceptional subgroup~\cite[Introduction]{mazu:mcei}.  The following
theorem now follows from Propositions~\ref{pr:borel},\ref{pr:split},
and \ref{pr:ns}.

\begin{thm}
Let $K$ be an imaginary quadratic field and $d$ a square-free positive integer.
There exists an effective constant $M_{K,d}$ such that, for all primes $p >
M_{K,d}$ and all $\Q$-curves $E/K$ of degree $d$, either
\begin{itemize} 
\item
the representation
\beq
\P\bar{\rho}_{E,p}: \GalQ \ra \PGL_2(\F_p)
\eeq
is surjective, or
\item $E$ has potentially good reduction at all primes not dividing
$6$.
\end{itemize}
\end{thm}

\section{Solutions to $A^4 + B^2 = C^p$}

We now use the results of the previous sections to study solutions to
the generalized Fermat equation
\begin{equation}
A^4 + B^2 = C^p
\label{eq:a4b2cp}
\end{equation}

The goal of this section is to prove the following theorem.

\begin{thm}
Suppose $A,B,C$ are coprime integers such that 
\beq
A^4 + B^2 = C^p
\eeq
and $p \geq 211$.  Then $AB = 0$.
\label{th:main}
\end{thm}

Suppose $(A,B,C)$ is a solution to \eqref{eq:a4b2cp} which is {\em
primitive} (i.e., $(A,B) = 1$) and which is {\em non-trivial} (i.e.,
$AB \neq 0$.)  We associate to $(A,B,C)$ a curve $E = E_{A,B,C}/\Q[i]$ with the
Weierstrass equation
\begin{equation}
E_{A,B,C}: y^2 = x^3 + 2(1+i)Ax^2 + (B+iA^2)x,
\label{eq:darmcurve}
\end{equation}
which was first discussed by Darmon in \cite{darm:sc} in connection with
the equation $A^4 + B^4 = C^p$.  We may think of $E$ as a
``generalized Frey-Hellegouarch curve'' whose relationship to
\eqref{eq:a4b2cp} is analogous to that between the usual
Frey-Hellegouarch curve and Fermat's equation.

Note that if $(A,B,C)$ is a solution to \eqref{eq:a4b2cp}, then so is
$(A,-B,C)$.  We therefore can and do assume that $B \equiv 0,2,3$ mod
$4$.

Write $\sigma$ for the non-trivial element of $\Gal(\Q[i]/\Q)$.  The map
\beq
\mu: (x,y) \mapsto
(\frac{1}{2}i(y^2/x^2),-\frac{1}{4}(1-i)y(B + iA^2-x^2)/x^2).
\eeq
is a degree $2$ isogeny from $E$ to its Galois conjugate $E^\sigma$.
Therefore, $E/\Q[i]$ is a $\Q$-curve of degree $2$.

One computes
\begin{eqnarray*}
E_4(E,\omega) & = & 80iA^2-48B \\
\Delta(E,\omega) & = & -64i(A^2+iB)(A^2-iB)^2,
\end{eqnarray*}
where $\omega$ is the Weierstrass differential $dx/2y$ with respect to
the Weierstrass equation \eqref{eq:darmcurve}.  Because $(A,B) = 1$,
we have that $E/\Q[i]$ is semistable away from $2$, and has
multiplicative reduction at an odd prime $\ic{p}$ of $\Q[i]$ precisely
when $\ic{p}$ divides $C$.

\begin{rem}  Suppose that $A$ and $B$ are chosen so that $A^4 + B^2$
is a prime number $\ell$.  Then $E$ has good reduction away from
primes of $\Q[i]$ dividing $2$ and $\ell$.  Moreover, the restriction
of scalars $\Res_{\Q[i]/\Q} E$ is an abelian surface over $\Q$ which
has good reduction away from $2$ and $\ell$.  We know, by the theorem of Iwaniec and
Friedlander~\cite{frie:a4b2}, that there are infinitely many
choices of $A,B$ such that $A^4 + B^2$ is prime; it follows that there
are infinitely many abelian surfaces over $\Q$ whose bad reduction is
supported at $2$ and a single odd prime.  It is interesting that we
know this fact for abelian surfaces, but not for elliptic curves over
$\Q$!
\end{rem}
 
We now embark on an analysis of the Galois representation
$\P\bar{\rho}_{E,p}$.  We will eventually show that, when $p$ is
large, this representation surjects onto $\PGL_2(\F_p)$.



We can define a lifting of $\P \bar{\rho}_{E,p}$ to an actual
representation as follows.   The abelian surface $A = \Res_{\Q[i]/\Q} E$ is an abelian surface with real multiplication by $\sqrt{2}$.  Let $\ic{p}|p$ be a prime of
$\Z[\sqrt{2}]$.  We define
\beq
\bar{\rho}_{E,\ic{p}}: \GalQ \ra \GL_2(\bar{\F}_p)
\eeq
to be the mod $p$ Galois representation attached to $A$.  Note that $\P
\bar{\rho}_{E,p}$ is, as the notation suggests, the projectivization
of $\bar{\rho}_{E,\ic{p}}$, and that
$\bar{\rho}_{E,\ic{p}}|\Gal(\bar{\Q}/\Q[i])$ is precisely the Galois representation $E[p](\bar{\Q})$.

We want to show, first of all, that $\P \bar{\rho}_{E,p}$ is
irreducible.  We begin with an elementary lemma on small primes.

\begin{lem} There exists a prime $\ell$ greater than $3$ which divides
$C$.
\label{le:smallprimes}
\end{lem}

\begin{proof}
If $A$ and $B$ were both odd, then $C^p$ would
be congruent to $2$ mod $4$, which is not possible.  So $C$ is odd.
Moreover, $C$ cannot be divisible by $3$, since $A^4 + B^2 = 0$ has no
nonzero solutions over $\F_3$.  Finally, $C \neq 1$, since $(A,B,C)$
is a non-trivial solution to \eqref{eq:a4b2cp}.  We conclude that
there exists a prime $\ell$ greater than $3$ which divides $C$; it
follows that $E$ has multiplicative reduction at primes of $\Q[i]$
over $\ell$.
\end{proof}

\begin{lem} $\bar{\rho}_{E,\ic{p}}$ is modular.
\label{le:emodular}
\end{lem}
\begin{proof} Since $3$ does not divide $C$, the curve $E$ has good reduction
at $3$, so $\P \bar{\rho}_{E,p}$ is unramified at $3$.  The modularity
of $E$ now follows from \cite[Th.\  5.2]{elle:qcurves}. 
\end{proof}

Proposition~\ref{pr:borel} and Lemma~\ref{le:smallprimes} imply that
$\bar{\rho}_{E,p}$ is irreducible.  Our next goal is to compute the
Serre invariants $N=N(\bar{\rho}_{E,\ic{p}}),
k=k(\bar{\rho}_{E,\ic{p}}),$ and $\epsilon =
\epsilon(\bar{\rho}_{E,\ic{p}})$.  By Ribet's theorem and Lemma~\ref{le:emodular}, we now have
\beq
\bar{\rho}_{E,\ic{p}} \cong \bar{\rho}_{f,p}
\eeq
for some $f$ in $S_k^{new}(N,\epsilon)$.  Note that $\det
\bar{\rho}_{E,\ic{p}}$ is cyclotomic, which implies that $\epsilon$ is
trivial and $k \cong 2$ mod $p-1$.

We now use the fact that $C^p = A^4 + B^2$ is a $p$th power.  This
fact implies that every odd prime $\ell$ dividing $\Delta(E,\omega)$
satisfies $p | \ord_\ell \Delta$.  By the theory
of the Tate curve, this implies that
$\bar{\rho}_{E,\ic{p}}|\Gal(\Qbar/\Q[i])$ is unramified away from $2$
and $p$, and so $\bar{\rho}_{E,\ic{p}}$ is unramified away from $2$
and $p$.  So $N$ is a power of $2$.

The fact that $p | \ord_\ic{q} \Delta$ for any prime $\ic{q} | p$
of $\Z[i]$ means, again using the Tate curve, that $E[p]/\Q[i]$
extends to a finite flat group scheme $\GG$ over the completion
$\Z[i]_\ic{q}$.  Since $\Z[i]/\Z$ is unramified at $p$, this extension
is unique.  By the \'{e}taleness of $\Z[i]/\Z$ at $p$, we can
descend $\GG \oplus \GG$ to a a finite flat group scheme over $\Z_p$
extending $A[p]/\Q$.  This means that 
$\bar{\rho}_{E,\ic{p}}$ is {\em finite} in the sense of
\cite[2.8]{serr:duke}, and so $k = 2$ by \cite[Prop. 4]{serr:duke}.

It remains to pin down $N$ precisely, which we accomplish by means of
Tate's algorithm.

\begin{prop} $N = 32$ if $A$ is even and $256$ if $A$ is odd.
\end{prop}

\begin{proof}
We begin by using Tate's algorithm to compute the local conductor $f$ of
the elliptic curve $E/\Q[i]$ at the prime $\pi = 1 + i$.  Recall that
$f$ is the multiplicity of the ideal $\pi$ in the conductor of $E$.
We refer to \cite[IV.\ \S 9]{silv:aec2} for all facts about Tate's algorithm.

First of all, we will translate $x$ by $1$, which gives us a new
Weierstrass equation
\beq
y^2 = x^3 + (3 + 2(1+i)A)x^2 + (3 + 4(1+i)A + B + iA^2) x
+ 2(1+i)A + 1 + B + iA^2.
\eeq
With respect to this equation, $b_2 = (8+8i)A + 12$.  If $A$ is odd
and $B$ even, then $\pi^2$ does not divide $a_6$, so $E$ has reduction
type II and $f = \ord_\pi(\Delta) = 12$.  (See \cite[IV.\ \S 9,Step
3]{silv:aec2}.)

If $A$ is even and $B$ odd, we define new variables by $y' = y - x$
and $x' = x - 1 - i$.  This change of variables gives rise to the
Weierstrass equation
\begin{eqnarray*}
(y')^2 + 2x'y'+ 2(1+i)y' = & &  \\
(x')^3 + (5 + 3i + 2(1+i)A)(x')^2 +
((B+iA^2) + (4 + 12i)A + 7 + 10i)x' \\ 
+ (B+iA^2)(i+2) + (-2 + 14i)A + 9i + 2. & & 
\end{eqnarray*}
 

We are now in the situation of \cite[IV.\ \S 9,Step 7]{silv:aec2}, so
$E$ has reduction of type $I^*_n$ for some $n$.  Note that $\pi^3 |
a_3$.  Also, $\pi^5$ divides $(-2 + 14i)A$ and $8i$, so we have $a_6
\cong (B+iA^2 + 1)(i+2)$ mod $\pi^5$.  Recall that we've assumed $B$
is not congruent to $1$ mod $4$.  So $B - 2A \equiv 3$ mod $4$.

If $B - 2A$ is congruent to $7$ mod $8$, we see
that $\pi^5 | a_6$. So the polynomial $Y^2 + \pi^{-2}a_3 Y - \pi^{-4}
a_6 Y$ has a double root over $\F_2$ at $Y = 0$.
Moreover, in this case
\beq
a_4 \cong B + 7 + 10i \cong B - 1 + 2i 
\eeq
mod $\pi^4$.  So $\ord_\pi a_4 = 3$, which implies that the polynomial
$\pi^{-1}a_2 X^2 + \pi^{-3} a_4 X + \pi^{-5} a_6$ has distinct roots
in $\F_2$.  We conclude in this case that $E$ has reduction type
$I_2^*$ and $f = \ord_\pi(\Delta) - 6 = 6$.

Suppose on the other hand that $B - 2A$ is congruent to $3$ mod $8$.
Then we change variables by setting $y'' = y' + 2$.  This change of
variables causes $a_6$ to become a multiple of $\pi^5$, while the
valuations of $a_4$ and $a_2$ do not change.  So once again we are in
the situation of reduction type $I_2^*$ and $f = 6$.

\medskip

The quantity $f$ computed above is the Artin conductor of $T_p E$
considered as a $p$-adic Galois representation of
$\Gal(\bar{\Q}_2/\Q_2[i])$.  Let $\rho_*$ be the the $4$-dimensional
representation of $\Gal(\bar{\Q}_2/\Q_2)$ induced from $T_p E$.  Then
$\rho_* \cong T_p A$, where $A$ is the restriction of scalars of $E$
described above.  We have from \cite[\S 1]{miln:arithav} that
$f(\rho_*) = f + 2 \dim T_p E = f + 4$.  So $f(\rho_*)$ is either $10$
or $16$.  

By examination of the $j$ invariant, we see that $E$ has potentially
good reduction at $\pi$.  It follows that the inertia group $I_2
\subset \Gal(\bar{\Q}_2/\Q_2)$ acts on $T_p A$ through a finite
quotient $G$, whose order is not divisible by any prime greater than
$3$.  Let $\rho_{E,\ic{p}}$ be the $2$-dimensional representation of
$\GalQ$ on $T_\ic{p} A$, and let $f_\ic{p}$ be the conductor of
$\rho_{E,\ic{p}} | \Gal(\bar{\Q}_2/\Q_2)$. If $\ic{p} = p$ is inert,
it is immediate that $f(\rho_*) = 2f_\ic{p}$.  If, on the other hand,
$\ic{p}$ and $\ic{p}'$ are split primes of $\Q(\sqrt{2})$ lying over
$p$, then $f_\ic{p}$ and $f_{\ic{p}'}$ are both equal to the $2$-part
of the conductor of the modular abelian variety $A$, and again we get
\beq
f(\rho_*) = f_\ic{p} + f_{\ic{p}'} = 2f_\ic{p}.
\eeq
We conclude that $f_\ic{p}$ is either $5$ or $8$.  Moreover, the fact
that $|G|$ is prime to $p$ implies that the $2$-part of the conductor
of $\bar{\rho}_{E,\ic{p}}$ is identical with $f_\ic{p}$.  This
completes the proof.
\end{proof}

\bigskip

We have now established that $\bar{\rho}_{E,\ic{p}}$ is isomorphic to
$\bar{\rho}_{f,p}$, where $f$ is a weight $2$ newform of level $32$ or
$256$.  In fact, the newforms of these levels are all associated to
elliptic curves (not necessarily defined over $\Q$) with complex
multiplication by $\Q[i]$ or $\Q[\sqrt{-2}]$.  In particular, the
image of $\P \bar{\rho}_{E,p}$ is the normalizer of a Cartan subgroup
in $PGL_2(\F_p)$.



Suppose the image of $\P \bar{\rho}_{E,p}$ lies in the normalizer
of a split Cartan subgroup.  Then it follows from
Proposition~\ref{pr:split} that $E$ has
good reduction away from $6$.  But this contradicts
Lemma~\ref{le:smallprimes}.

We conclude that the image of $\P \bar{\rho}_{E,p}$ must be
the normalizer of a non-split Cartan subgroup. We now use this fact
to bound $p$.

\begin{prop} Let $E_0/\Q[i]$ be a $\Q$-curve of degree $2$, and suppose
$\P \bar{\rho}_{E,p}$ has image contained in the normalizer of a
non-split Cartan subgroup of $\PGL_2(\F_p)$, for some $p \geq 211$.  Then $E_0$ has potentially good
reduction for all primes of $\Q[i]$.
\label{pr:211}
\end{prop}

Proposition~\ref{pr:211} completes the proof of
Theorem~\ref{th:main}.  For we have shown that if $(A,B,C)$ is a
solution to \eqref{eq:a4b2cp}, and $E = E_{A,B,C}$, then $\P
\bar{\rho}_{E,p}$ has as image the normalizer of a non-split Cartan
subgroup of $\PGL_2(\F_p)$.  If $p \geq 211$, then
Proposition~\ref{pr:211} shows that $E$ has
good reduction everywhere, contradicting
Lemma~\ref{le:smallprimes}.

We now proceed with the proof of Proposition~\ref{pr:211}.

\begin{proof}
The main tools are Proposition~\ref{pr:ns} and the estimate for
average special values of $L$-functions in \cite{elle:duke}.

Proposition~\ref{pr:ns} tells that $E_0$ has
potentially good reduction everywhere if $p$ is sufficiently
large.  We now show that $p \geq 211$ suffices.  It
is clear from the proof of Proposition~\ref{pr:ns} that $E_0$
has good reduction away from $6p$ whenever there exists a newform $f$
on $S_2(\Gamma_0(p^2))$ satisfying the conditions of
Proposition~\ref{pr:analytic}.  In turn, in order to prove the
existence of such a form it suffices to show that the inner product
\beq
(a_1, L_\chi)_{p^2}^{p-new} = (a_1,L_\chi)_{p^2} - (a_1, L_\chi)^{p}_{p^2}
\eeq
is non-zero. We have by Lemma~\ref{le:fromdp} that 
\beq
(a_1, L_\chi)^{p}_{p^2} = 
p(p^2-1)^{-1}(a_1 - p^{-1}\chi(p)a_p, L_\chi)_p.
\eeq

We now use Theorem 1 of \cite{elle:duke} to show that
$(a_1, L_\chi)^{p-new}$ is nonzero.  The theorem shows that
$(a_m,L_\chi)_N$ is $4\pi\chi(m) + O(N^{-1}\log(N)d(N))$, where
$d(N)$ is the number of divisors of $N$ andthe
implied constants are explicit functions of $m,\chi$.  Precisely, we
obtain that, for $p \geq 211$, 
\beq
|(a_1,L_\chi)_{p^2} - 4\pi| \leq 4.37
\eeq

(In the applications here, we always take the parameter $\sigma$ in
\cite{elle:duke} to be $8/\pi$.)

At level $p$, the same theorem gives
\beq
|(a_1,L_\chi)_{p}| \leq 786
\eeq
when $p > 211$.

Finally, Lemma~\ref{le:convexity} shows that
\beq
|(a_p,L_\chi)_{p}| \leq 437.
\eeq
So we find
\beq
|(a_1,L_\chi)^{p-new}| \geq 4\pi - 4.37 - \frac{211}{211^2-1}(786 +
437/211) > 4.
\eeq
This proves Proposition~\ref{pr:211}.
\end{proof}

For $p$ smaller than $211$, the argument above shows that $A^4 + B^2 = C^p$
has no nontrivial solutions if we can prove the existence of a modular
form satisfying the conditions of Proposition~\ref{pr:analytic}.  To
be precise, we have shown

\begin{prop}  Let $p>13$ be prime, and suppose there exists either
\begin{itemize}
\item a newform in $S_2(\Gamma_0(2p^2))$ with $w_p f = f$ and $w_2 f =
-f$; or 
\item a newform in $S_2(\Gamma_0(p^2))$ with $w_p f = f$,
\end{itemize}
such that $L(f \tensor \chi, 1) \neq 0$. Then the equation $A^4 + B^2
= C^p$ has no primitive non-trivial solutions.
\label{pr:finitecomp}
\end{prop}

Verification of the existence of such a modular form is, in principle,
a finite computation.  In practice, it is beyond the reach of current
computers to compute the Fourier coefficients of a newform of level
$2p^2$ when $p$ is as large as $100$.  It seems probable that by
exploiting various tricks and carrying out more complicated
computations, we will be able to show that
Proposition~\ref{pr:finitecomp} applies for all primes $p$ between $17$
and $211$.  We will discuss this problem in a later paper.

\end{document}